\documentclass[11pt]{article}

\usepackage{amsmath,amssymb}

\newtheorem{theorem}{Theorem}
\newtheorem{ass}{Theorem}[section]
\newtheorem{prop}[ass]{Proposition}
\newtheorem{lemma}[ass]{Lemma}

\newtheorem{prob}{Problem}

\newcommand{\qed}{\hspace*{\fill} \rule{7pt}{7pt}}
\newcommand{\Proof}{\noindent{\bf Proof.}\ \ }
\newcommand{\C}{{\cal C}}

\begin{document}

\title{Packing without some pieces}

\author{
Raphael Yuster
\thanks{Department of Mathematics, University of Haifa, Haifa
31905, Israel. Email: raphy@math.haifa.ac.il.
This research was supported by the Israel Science Foundation
(grant No. 1082/16).}
}

\date{}

\maketitle

\setcounter{page}{1}

\begin{abstract}

Erd\H{o}s and Hanani proved that for every fixed integer $k \ge 2$, the complete graph $K_n$ can be almost completely packed 
with copies of $K_k$; that is, $K_n$ contains pairwise edge-disjoint copies of $K_k$ that cover all but an $o_n(1)$ fraction of
its edges. Equivalently, elements of the set $\C(k)$ of all red-blue edge colorings of $K_k$ can be used to
almost completely pack every red-blue edge coloring of $K_n$.

The following strengthening of the aforementioned Erd\H{o}s-Hanani result is considered.
Suppose $\C' \subset \C(k)$. Is it true that we can use elements only from $\C'$
and almost completely pack every red-blue edge coloring of $K_n$?
An element $C \in \C(k)$ is {\em avoidable} if $\C'=\C(k) \setminus C$ has this property
and a subset ${\cal F} \subset \C(k)$ is avoidable if $\C'=\C(k) \setminus {\cal F}$ has this property.

It seems difficult to determine all avoidable graphs as well as all avoidable families.
We prove some nontrivial sufficient conditions for avoidability. Our proofs imply, in particular, that
(i) almost all elements of $\C(k)$ are avoidable
(ii) all Eulerian elements of $\C(k)$ are avoidable
and, in fact, the set of all Eulerian elements of $\C(k)$ is avoidable.

\vspace*{5pt}

\end{abstract}

\section{Introduction}

Throughout this paper a red-blue edge coloring of $K_k$ is synonymous with a graph $H$ on $k$ vertices where $E(H)$ are the blue
edges and $E(H^c)$ are the red edges. We usually omit the word ``edge'' and just refer to red-blue colorings.
Let $\C(k)$ be the set of all red-blue colorings of $K_k$. Equivalently, we can view $\C(k)$ as the set of all graphs on $k$ vertices.

If $F_1, F_2, \ldots, F_t$ are pairwise edge-disjoint cliques of size $k$ forming a packing of $K_n$, then given any red-blue coloring of $K_n$ with color classes
$G_{blue}$ and $G_{red}$, we can view the $F_i$'s as red-blue colorings of $K_k$ where the coloring of $F_i$ is given by $F_i \cap G_{blue}$ and
$F_i \cap  G_{red}$ for $i = 1,\ldots,t$. The main question of the paper is what possible $2$-colorings $F_i \cap G_{blue}$  and $F_i \cap  G_{red}$ are forced to arise in asymptotic packings (packings that cover almost all of the edges of $K_n$).

More formally,
for $X \subseteq \C(k)$ an {\em $X$-packing} of a red-blue coloring of $K_n$ is a set ${\cal P}$ of pairwise edge-disjoint
subgraphs of this colored $K_n$, where each subgraph is isomorphic to an element of $X$. The size of the packing is $|{\cal P}|$.
Obviously, $|{\cal P}| \le \frac{n(n-1)}{k(k-1)}$.

We say that $X$ has the {\em asymptotic packing property} if every red-blue coloring of $K_n$ has an
$X$-packing of size at least $\frac{n(n-1)}{k(k-1)}(1-o_n(1))$. More formally, for every $\epsilon > 0$
and all sufficiently large $n$,
there is an $X$ packing of every red-blue coloring of $K_n$ of size at least $\frac{n(n-1)}{k(k-1)}(1-\epsilon)$.
The following was proved by Erd\H{o}s and Hanani \cite{EH-1963}:
\begin{theorem}\label{t:eh}
\C(k) has the asymptotic packing property.
\end{theorem}
In other words, they proved that $K_n$ can be packed with edge-disjoint copies of $K_k$ so that only $o(n^2)$ edges remain unpacked.
This result has many applications and was generalized in several ways, most notably by R\"odl for hypergraphs \cite{rodl-1985},
by Wilson for exact graph decompositions \cite{wilson-1975} and by Keevash for exact hypergraph decompositions \cite{keevash-2014}.
See also Glock et al. \cite{GKLO-2016} for another, more general proof.

It is therefore interesting to determine to what extent can Theorem \ref{t:eh} be strengthened by requiring less than $\C(k)$
in its statement.
Namely, which subsets of $\C(k)$ have the asymptotic packing property.
\begin{prob}\label{prob:1}
For every fixed $k$, determine the subsets of $\C(k)$ that have the asymptotic packing property.
\end{prob}
An element $C \in \C(k)$ is {\em avoidable} if $\C'=\C(k) \setminus C$ has the asymptotic packing property
and a subset ${\cal F} \subset \C(k)$ is avoidable if $\C'=\C(k) \setminus {\cal F}$ has the asymptotic packing property.
Non-avoidable graphs or subsets are {\em unavoidable}. So Problem \ref{prob:1} can be reformulated
as asking to determine all avoidable subsets and in particular all avoidable graphs.

For $k=2$ we trivially have that every nonempty subset of $\C(2)$ is unavoidable.
It is also easy to verify that every nonempty subset of $\C(3)$ is unavoidable. In fact:
\begin{prop}\label{p:1}
For all $k \ge 2$, the graphs $K_k$, $K_{1,k-1}$ and their complements are unavoidable.
Also, $K_{2,3}$, $K_{3,4}$ and $K_4^{-}$ and their complements are unavoidable.
\end{prop}
Already for $k=4$ we do not know the complete solution for Problem \ref{prob:1}.

Let ${\cal U}(k) \subseteq \C(k)$ denote the set of all unavoidable graphs on $k$ vertices.
Our first main result is that almost all elements of $\C(k)$ are avoidable.
\begin{theorem}\label{t:1}
$|{\cal U}(k)| = o(|\C(k)|)$.
\end{theorem}
Theorem \ref{t:1} is a consequence of a result that gives a more general sufficient condition for avoidability
in terms of the asymmetry of a graph (Lemma \ref{l:1}). It is natural to use random $k$-vertex graphs as it is not difficult to prove that these are
almost surely highly asymmetric (in a well-defined sense made later). The main technical issue is proving that this asymmetry property suffices for
avoidability.

While Theorem \ref{t:1} shows that graphs that are sufficiently asymmetric are avoidable, our second main result proves
that a certain large class of graphs which contains some highly symmetric graphs is avoidable. This class of graphs, whose definition follows,
includes all Eulerian elements of $\C(k)$.

The {\em degree set} of a graph $G$ is the set $\{d(v) ~|~ v \in V(G)\}$.
For a set of integers $S \subseteq \{0,\ldots,k-1\}$ let ${\cal F}(S,k)$ be the set of all graphs on $k$ vertices whose degree
set is contained in $S$. So, ${\cal F}(\{t\},k)$ is the set of all $t$-regular graphs on $k$ vertices.
Equivalently, ${\cal F}(S,k)$ is the set of all red-blue colorings of $K_k$ where the degree set of each blue graph is contained
in $S$. When $k$ is odd, a red-blue coloring of $K_k$ is Eulerian if the blue graph is Eulerian and the red graph is
Eulerian. For example, a coloring of $K_5$ with a blue $C_5$ (and hence a red $C_5$) is Eulerian.
Notice that all Eulerian red-blue colorings are contained in ${\cal F}(S,k)$ where $S=\{2,4,\ldots,k-3\}$, but the latter is
more general already for $k=7$.
An immediate corollary of the following theorem is that the family of all Eulerian red-blue edge-colorings of $K_k$ is avoidable.

\begin{theorem}\label{t:2}
For all odd positive integers $k$, ${\cal F}(\{2,4,\ldots,k-3\},k)$ is avoidable.
\end{theorem}
Theorem \ref{t:2} is a nontrivial consequence of a more general statement (Theorem \ref{t:suff}) that gives a sufficient condition for the
avoidability of ${\cal F}(S,k)$ in terms of the solvability of a certain parametric linear program.
For relatively small $k$ we can determine if a solution exists and hence determine many additional $S$ such that ${\cal F}(S,k)$
is avoidable.

The tool of fractional packings will be useful in proving Theorem \ref{t:1}, Theorem \ref{t:2}, and their more generalized statements.
We describe this tool in Section 2. Sections 3 and 4 prove Theorem \ref{t:1} and Theorem \ref{t:2} respectively.
Section 5 contains the proof of Proposition \ref{p:1}. The final section contains some concluding remarks,
most notably  addressing the analogous problem where instead of an asymptotic packing we ask for an {\em exact} decomposition
and consider the seemingly stronger property of decomposition avoidability. In particular, we prove there that $C_4$ is not decomposition avoidable.

\section{Fractional packings}

Let ${\cal R}$ be a set of graphs of order $k$. Let $G$ be a graph with $V(G)=[n]$. Let $\binom {G}{\cal R}$ denote the set of all induced copies
of ${\cal R}$ in a graph $G$ (by induced copy we mean an induced subgraph of $G$ on $k$ vertices which is isomorphic to an element of ${\cal R}$).
Notice that in the special case that ${\cal R}$ contains all induced $k$-subgraphs of $G$,
then $|\binom {G}{\cal R}|=\binom{n}{k}$. 

A function $\phi$ from $\binom {G}{\cal R}$ to $[0,1]$ is a {\em fractional ${\cal R}$-packing} of $G$ if for each pair of distinct vertices
$\{x,y\} \subset [n]$ we have
\begin{equation}\label{e:fp}
\sum_{H \in \binom {G}{\cal R} \,:\, \{x,y\} \subset V(H)} {\phi(H)} \leq 1\;.
\end{equation}
For a fractional ${\cal R}$-packing $\phi$, let
$$
|\phi|=\sum_{H \in \binom{G}{\cal R}} \phi(H)\;.
$$
The {\em fractional ${\cal R}$-packing number}, denoted by $\nu^*_{\cal R}(G)$, is the maximum value
of $|\phi|$ ranging over all fractional ${\cal R}$-packings $\phi$.
One observes that computing $\nu^*_{\cal R}(G)$ amounts to solving a linear programming maximization problem with
$\binom{n}{2}+|\binom {G}{\cal R}|$ constraints and $|\binom {G}{\cal R}|$ variables.
It can therefore be solved in polynomial time for fixed $k$.

An ${\cal R}$-packing of $G$ is a {\em fractional ${\cal R}$-packing} whose image is $\{0,1\}$.
In other words, it is a set of induced copies of elements of ${\cal R}$ in $G$ where any two copies do not share a pair of vertices
(they are either disjoint or have a single vertex in common).
Let $\nu_{\cal R}(G)$ denote the maximum size of an ${\cal R}$-packing of $G$.
As we restrict the values of $\phi$ in the definition of an ${\cal R}$-packing of $G$, we have $\nu^*_{\cal R}(G) \geq \nu_{R}(G)$.

An important result of Haxell and R\"odl \cite{HR-2001} and later a slightly more general form (allowing for a ``set of graphs'' definition) by
the author \cite{yuster-2005}, both of which rely on Szemer\'edi's regularity lemma \cite{szemeredi-1978}, shows that the converse
inequality is also asymptotically true, up to an additive error term which is negligible for dense graphs.

\begin{lemma}\label{l:hr}
For every $\epsilon > 0$ and for every positive integer $k \ge 2$ there exists $N=N(k,\epsilon)$ such that the following holds.
For any set ${\cal R}$ of graphs of order $k$ and any graph $G$ with $n > N$ vertices, $\nu^*_{\cal R}(G)-\nu_{\cal R}(G) \le \epsilon n^2$.
\end{lemma}

One can observe that Lemma \ref{l:hr} is extremely useful already by the following trivial use of it which implies the (nontrivial) result of
Erd\H{o}s and Hanani. Indeed, merely notice that if ${\cal R}=\{K_k\}$ and $G=K_n$, then clearly $\nu^*_{\cal R}(K_n)=\binom{n}{2}/\binom{k}{2}$.
Thus, $\nu_{\cal R}(G)=\binom{n}{2}/\binom{k}{2}-o(n^2)$.

\section{Avoidable graphs}

\subsection{Decompositions and fractional decompositions}

We say that $X \subseteq \C(k)$ has the {\em decomposition property for $n$} if every red-blue coloring of $K_n$ has an
$X$-packing of size $\frac{n(n-1)}{k(k-1)}$. Notice that having the decomposition property for $n$ is the same as having
$\nu_{X}(G)=\frac{n(n-1)}{k(k-1)}$ for every graph $G$ with $n$ vertices. Analogously, we say that $X$ has the
{\em fractional decomposition property for $n$} if $\nu^*_{X}(G)=\frac{n(n-1)}{k(k-1)}$. Trivially, $\C(k)$ has the fractional decomposition property
for all $n \ge k$, and a seminal result of Wilson \cite{wilson-1975} asserts that $\C(k)$ has the decomposition property
for all $n$ sufficiently large that satisfy the necessary divisibility condition $n \equiv 1,k \bmod k(k-1)$.

Let $H$ be a graph with $h$ vertices. For $1 \le k \le h$, let $\C(H,k)$ be the set of all induced subgraphs of $H$ on $k$ vertices.
So, for example, if $H=C_6$ and $k=4$, then $\C(C_6,4)=\{P_4,P_3 \cup K_1,2K_2\}$.
\begin{lemma}\label{l:reduce}
Let $H$ be a graph with $h$ vertices. Suppose that $X=\C(k) \setminus \C(H,k)$ has the decomposition property for some $q$.
Then $H$ is avoidable.
\end{lemma}
\Proof
Let $k \le h$ be maximal such that $X=\C(k) \setminus \C(H,k)$ has the decomposition property for some $q$.
Let $q$ be minimal subject to this, so $q=q(H)$ only depends on $H$.

Consider first the easy case where $k=h$. In this case already $X=\C(h) \setminus H$ has the decomposition property for $q$.
Then we can decompose every red-blue coloring of $K_q$ into pairwise edge-disjoint copies of $K_h$ where in each copy,
the blue edges do not induce $H$. By Theorem \ref{t:eh} (the Erd\H{o}s-Hanani Theorem), $C(q)$ has the asymptotic packing
property. Thus, $K_n$ can be packed with edge-disjoint copies of $K_q$ so that only $o(n^2)$ edges remain unpacked.
This, in turn, implies that any red-blue coloring of $K_n$ can be packed with edge-disjoint copies of $K_h$ so that only $o(n^2)$ edges remain
unpacked, and in each copy, the blue edges do not induce $H$.  Thus, $\C(h) \setminus H$ has the asymptotic packing property, which means
that $H$ is avoidable.

Now consider the case where $k < h$. By the result of Wilson mentioned earlier, there exists $n_0=n_0(q)=n_0(H)$ such that for all
$n > n_0$, if $n \equiv 1 \bmod q(q-1)$, then $K_n$ has a decomposition into $\frac{n(n-1)}{q(q-1)}$ pairwise edge-disjoint copies of $K_q$,
and as in the previous paragraph, each such $K_q$ can be decomposed into $\frac{q(q-1)}{k(k-1)}$ pairwise edge-disjoint copies of $K_k$,
such that in each copy of $K_k$, the blue edges are isomorphic to an element of $X$. Altogether, any red-blue coloring of $K_n$
has a decomposition into $\frac{n(n-1)}{k(k-1)}$ pairwise edge-disjoint copies of $K_k$, such that in each copy of $K_k$, the blue edges are
isomorphic to an element of $X$. Let ${\cal D}$ denote the elements of this decomposition.

But recall that we want to prove that $H$ is avoidable (and not merely that $\C(H,k)$ is avoidable).
To this end, let us design some fractional packing of $K_n$. Consider some $K \in {\cal D}$ and recall that $K \in X$.
There are $n-k$ vertices of $K_n$ that do not belong to $K$. For any set $T$ of $h-k$ of these vertices (there are $\binom{n-k}{h-k}$ choices for
$T$) consider the $K_h$-subgraph of $K_n$ induced by the vertices of $K$ and the vertices of $T$, call it $Y$.
Notice that $Y$ is a red-blue coloring of $K_h$ where the blue edges of $Y$ do not induce a subgraph that is
isomorphic to $H$. Indeed, this is because $K$ is an induced $k$-vertex subgraph of $Y$, so if $Y$ were isomorphic to $H$,
then $K$ would have been a member of $\C(H,k)$ while by definition $K \in X=\C(k) \setminus \C(H,k)$.
We give $Y$ the weight $x$ ($x$ to be chosen later). We do this for every choice of $K \in {\cal D}$ and for every choice of $T$,
and they all get the same weight $x$. So, altogether we obtain a fractional packing of $K_n$ consisting of
$$
\frac{n(n-1)}{k(k-1)} \cdot \binom{n-k}{h-k}
$$
elements, each one  having weight $x$, and each one being a red-blue coloring of $K_h$ with the blue edges not forming an $H$.
Since, by symmetry, the sum of the weights of each edge of $K_n$ is the same, we can choose the weight $x$ such that
the total weight of this fractional packing is precisely $\frac{n(n-1)}{h(h-1)}$ (a fractional decomposition).
In other words, $\nu^*_R(G) = \frac{n(n-1)}{h(h-1)}$ where $R=C(h) \setminus H$ and $G$ is any graph on $n$ vertices.

There are still two small issues to take care of. First observe that the argument above assumed that $n \equiv 1 \bmod q(q-1)$
(and recall that $q=q(H)$). If $n > n_0$ is not of this form, that let $n' < n$ be the largest integer such that $n' \equiv 1 \bmod q(q-1)$.
As $n-n' \le q(q-1) = o(n)$, we can just ignore $n-n'$ vertices (which touch $o(n^2)$ edges) and thus
$\nu^*_R(G) = \frac{n(n-1)}{h(h-1)}-o(n^2)$ where $R=C(h) \setminus H$ and $G$ is any graph on $n$ vertices.
Finally, we can use Lemma \ref{l:hr} to obtain that $\nu_R(G) = \frac{n(n-1)}{h(h-1)}-o(n^2)$.
But this means that $H$ is avoidable, as required.
\qed

In the proof of Theorem \ref{t:1} it would be very important to use Lemma \ref{l:reduce} for $k$ which is very close to $h$ and for $q$
which is not too large. Quantitatively, this will be guaranteed by the following lemma.
\begin{lemma}\label{l:dependence}
For every $h \ge 2$, there is $h-o(h) \le k \le h$ such that $K_q$ has a decomposition into $q$
pairwise edge-disjoint copies of $K_k$ where $q=k^2-k+1$.
\end{lemma}
\Proof
Assume first that $p=h-1$ is a prime power. It is well-known that there is a finite projective plane of order $p$, which means
that $K_{p^2+p+1}$ decomposes into $p^2+p+1$ pairwise edge-disjoint copies of $K_h$. So, in this case, the lemma holds for $k=h$.
In the case where $h-1$ is not a prime power, we can use the result of Baker, Harman, and Pintz \cite{BHP-2001} which 
states that there is always a prime strictly between $x$ and $x+O(x^{21/40})=x+o(x)$ (this result is a significant extension of
Chebyshev's Theorem of Bertrand's  postulate).
So, let $k \le h$ be the largest integer such that $k-1$ is a prime power.
Since $k \ge h-o(h)$, using the same argument of existence of projective plane of order $k-1$, we have that $K_q$ decomposes into
$q$ pairwise edge-disjoint copies of $K_k$ where $q=k^2-k+1$.
\qed

It is important to note that if we wouldn't have cared about the fact that $q$ is small (only a polynomial in $h$),
then Lemma \ref{l:dependence} would have worked already with $k=h$ since Wilson's Theorem mentioned earlier guarantees that for
some large $q$, $K_q$ has a decomposition into $K_h$. However, the bound in Wilson's proof for such a $q$ does not suffice for our proof.

\subsection{Graphs whose large subgraphs are asymmetric}

Let $H$ be a graph on the vertex set $[h]=\{1,\ldots,h\}$.
A permutation $\pi : [h] \rightarrow [h]$ is an {\em automorphism} of $H$
if $(\pi(i), \pi(j))$ is an edge of $H$ if and only if $(i,j)$ is an edge of $H$.
The group of all automorphisms of $H$ is denoted by $aut(H)$.
We say that $H$ is {\em asymmetric} if $aut(H)$ consists only of the identity permutation.
Otherwise, we say that $H$ is {\em symmetric}.
The smallest graph (with more than one vertex) which is asymmetric
is obtained from the path on vertices $1,2,3,4,5$ (in this order) by adding vertex $6$
and connecting it to vertices $3$ and $4$.
Erd{\H{o}}s and R\'enyi \cite{ER-1963} proved that almost all graphs are asymmetric.

For a graph $H$, let $k(H)$ be the smallest integer $k$ such that
any two induced subgraphs of $H$ on at least $k$ vertices each, are non-isomorphic
and further, any induced subgraph of $H$ on at least $k$ vertices is asymmetric.
If $H$ is symmetric then define $k(H)=\infty$.

It is not difficult to prove that $k(H) \ge \lceil (h+1)/2 \rceil$ as it is well-known (Goodman \cite{goodman-1959})
that for any graph, there are two vertices that {\em agree} on at least
$\lceil (h-3)/2 \rceil$ other vertices, where $u$ and $v$ agree on $w$ if both are neighbors of $w$
or both are non-neighbors of $w$.
Asymmetric graphs are natural candidates for a graph with relatively small $k(H)$, but this is
clearly not a sufficient condition, as it is easy to construct asymmetric graphs with $k(H)=h-o(h)$.
We will need graphs $H$ with relatively small $k(H)$ as it would be possible to
prove that such $H$ are avoidable.

Our next lemma proves that a randomly chosen graph $H$ on $h$ vertices has relatively
small $k(H)$, with probability tending to one as $h$ increases.
Recall that ${\cal G}(h,\frac{1}{2})$ is the probability space of all graphs on $h$
vertices where each pair of vertices are connected with an edge with probability
$\frac{1}{2}$, and the $\binom{h}{2}$ choices are independent.

\begin{lemma}\label{l:low-k}
Let $\beta \ge 0.94$ be fixed and let $H \sim {\cal G}(h,\frac{1}{2})$. Then,
$$
\Pr \left[k(H) \le \beta h\right] = 1-o_h(1)\;.
$$
\end{lemma}
\Proof
Recall that $V(H)=[h]$. For a subset $K \subseteq [h]$ let $H[K]$ be the subgraph of $H$ induced by $K$.
We will prove the following two claims.
\begin{enumerate}
\item[C1.]
For every $K \subseteq [h]$ with $|K| \ge \beta h$, the probability that $K[H]$ is symmetric is at most $(1.3)^{-h}$.
\item[C2.]
For any two distinct subsets $J,K \subseteq [h]$ with $|J|=|K| \ge \beta h$ such that $|J \triangle K| = 2t$, the probability that $K[H]$ and $J[H]$ are isomorphic is at most $(1.3)^{-ht}$.
\end{enumerate}
There are less than $h \binom{h}{\beta h}$ subsets $K \subseteq [h]$ of size at least $\beta h$.
For any such $K$, the number of subsets $J$ with $|J|=|K|$ such that $|J \triangle K| = 2t$ is less
than $h^{2t}$. Also notice that $ 2 \le |J \triangle K| \le 2h-2|K| \le 0.12h$ so $t \le 0.06h$.
Thus, if both claims hold we obtain by the union bound that
$$
\Pr \left[k(H) \le \beta h\right] \ge 1- h \binom{h}{\beta h} (1.3)^{-h} - h \binom{h}{\beta h} \left(\sum_{t=1}^{\lceil 0.06 h \rceil} h^{2t} (1.3)^{-ht}\right) \ge 1-o_h(1)
$$
where in the last inequality we have used the fact that $\beta \ge 0.94$ which implies that
$\binom{h}{\beta h} = o((1.26)^h)$.

We next prove Claim C1.
Let $K \subseteq [h]$ with $|K|=k \ge \beta h$.
Let $\pi$ be a permutation of $K$ which is not the identity.
We would like to upper bound the probability that $\pi \in aut(H[K])$.
As in the proof of Kim, Sudakov, and Vu \cite{KSV-2002}, it would be useful to compute such a bound by considering
the number of non-stationary points of $\pi$. Let this number be $s$.
Notice that for any given $s$, the number of possible $\pi$ with $s$ non-stationary points is less than $k^s$.

Let $S \subseteq K$ be the set of non-stationary points of $\pi$, so $|S|=s \ge 2$.
Observe that we can always find $r = \lceil s/3 \rceil$ pairs $\{v_1,u_1\},\ldots,\{v_r,u_r\}$ such that
$v_1,\ldots,v_r,u_1,\ldots,u_r$ are distinct elements of $S$ and further $\pi(v_i)=u_i$ for $i=1,\ldots,r$.
Indeed, in each nontrivial orbit of $\pi$ of length $\ell$ we can  obviously find $\lfloor \ell/2 \rfloor$ such pairs.
The worst case is when all nontrivial orbits are of length $3$ so we can only find one pair in each orbit,
resulting in only $s/3$ pairs. Let $S^*=\{v_1,\ldots,v_r,u_1,\ldots,u_r\} \subseteq S$.

For each $i=1,\ldots,r$ and for each point $w \in K \setminus S^*$ consider the two pairs
$\{v_i,w\}$ and $\{u_i,\pi(w)\}$. Since $\pi(v_i)=u_i$, in order for $\pi$ to
be in $aut(H[K])$ we must have that $\{v_i,w\}$ and $\{u_i,\pi(w)\}$ {\em agree} (both are edges or both are non-edges).
Since agreement occurs with probability $\frac{1}{2}$ and since all the $r \cdot (k-2r)$ choices of $i$ and $w$
are independent with respect to the event of agreement (since they correspond to distinct pairs), we obtain that
\begin{equation}\label{e:pi}
\Pr[\pi \in aut(H[K]) ] \le 2^{-r(k-2r)} = 2^{-\lceil s/3 \rceil (k-2\lceil s/3 \rceil)}\;.
\end{equation}
Now, since $k \ge \beta h \ge 0.94 h$, we obtain that for all $s=2,\ldots,k$,
$$
2^{-\lceil s/3 \rceil (k-2\lceil s/3 \rceil)} k^s < \frac{1}{k (1.3)^h}\;.
$$
Notice that for $h$ sufficiently large, the left hand side is maximized when $s=3$
and already in this case the inequality holds since $2^{0.94} > 1.3$.

As there are less than $k^s$ permutations $\pi$ with $s$ non-stationary points we obtain by the union bound,
the last inequality, and (\ref{e:pi}) that
$$
\Pr[H[K] \text{ is symmetric} ] \le \sum_{s=2}^k k^s 2^{-\lceil s/3 \rceil (k-2\lceil s/3 \rceil)} < \sum_{s=2}^k \frac{1}{k(1.3)^h} <
\frac{1}{(1.3)^h}\;.
$$
This completes the proof of Claim C1.

We next prove Claim C2 which is quite similar.
Let $J,K \subseteq [h]$ with $|J|=|K| \ge \beta h$ such that $|J \triangle K| = 2t$.
Let $\pi$ be a bijection from $K$ to $J$.
We would like to upper bound the probability that $\pi$ is an isomorphism between $H[K]$ and $H[J]$.
Let $s$ be the number of non-stationary points of $\pi$. Observe that $s \ge |K \setminus J|=t$
and that the number of possible $\pi$ with $s$ non-stationary points is at most $k^s$.

We claim that we can always find $r=\min\{\lfloor k/4 \rfloor, t+ \lceil (s-t)/3 \rceil \}$ pairs
$\{v_1,u_1\},\ldots,\{v_r,u_r\}$ such that all the $2r$ vertices are distinct,
$v_i \in K$, $u_i \in J$ and $\pi(v_i)=u_i$.
Indeed, the vertices of $K \setminus J$ are all non-stationary, so we let them be $v_1,\ldots,v_t$
and let their images be $u_1,\ldots,u_t$, respectively. Each $u_i$ may be either in $J \setminus K$
or in $J \cap K$. suppose $m$ of them are in $J \setminus K$. Then there are $t-m$ additional
vertices $v_{t+1},\ldots,v_{2t-m}$ in $J \cap K$ having images in $J \setminus K$.
Denote the images by $u_{t+1},\ldots,u_{2t-m}$ respectively. This still leaves
 $s-2t+m$ non-stationary vertices of $J \cap K$ having images also in $J \cap K$, so as in the proof of Claim
C1, we can pick at least $\lceil (s-2t+m)/3 \rceil$ additional pairs $\{v_i,u_i\}$ such that $\pi(v_i)=u_i$
for $i=2t-m+1,\ldots,2t-m+ \lceil (s-2t+m)/3 \rceil$ and such they are all distinct from the vertices
in previously selected pairs. So, the least amount of selected pairs occurs when $m=t$
in which case we can still pick at least $t+\lceil (s-t)/3 \rceil$ pairs.
This proves the claim about the existence of $r$. The reason we take $r$ to be the minimum between
$t+ \lceil (s-t)/3 \rceil$ and $\lfloor k/4 \rfloor$ is that we still want to leave sufficiently many vertices
of $K$ that are not in these $r$ pairs.
Let $S^*=\{v_1,\ldots,v_r,u_1,\ldots,u_r\}$ and observe that $|K \setminus S^*| \ge k-2r$.

For each $i=1,\ldots,r$ and for each point $w \in K \setminus S^*$ consider the two pairs
$\{v_i,w\}$ and $\{u_i,\pi(w)\}$. Since $\pi(v_i)=u_i$, in order for $\pi$ to
be an isomorphism we must have that $\{v_i,w\}$ and $\{u_i,\pi(w)\}$ agree (both are edges or both are non-edges).
Since agreement occurs with probability $\frac{1}{2}$ and since all the $r |K \setminus S^*| \ge r(k-2r)$ choices of $i$ and $w$
result in independent events (since they correspond to distinct pairs), we obtain that
\begin{equation}\label{e:pi-2}
\Pr[\pi \text{ is an isomorphism}] \le 2^{-r(k-2r)} \le 2^{-rk/2}
\end{equation}
where in the last inequality we have used the fact that $r \le k/4$.

Consider first the case where $r=\lfloor k/4 \rfloor$.
In this case we have for all $h$ sufficiently large that
$$
2^{-rk/2} k^s \le 2^{-k^2/9}k^s \le 2^{-(0.94h)^2/9}k^h < \frac{1}{k(1.3)^{0.06h^2}} < \frac{1}{k(1.3)^{ht}}
$$
where we have used here the fact that $k \ge 0.94h$ and $t \le 0.06h$.

Consider next the remaining case where $r=t+ \lceil (s-t)/3 \rceil \} \le \lfloor k/4 \rfloor$.
Let $b = \lceil (s-t)/3 \rceil$ so $r=t+b$ and $s \le 3b+t$.
We now have that
\begin{eqnarray*}
2^{-rk/2} k^s & = & 2^{-(t+ b)k/2}k^s \\
& \le & 2^{-(t+b)0.47h}k^{3b+t} \\
& =  & 2^{-0.47ht} 2^{-b(0.47h-3\log_2 k)}2^{t \log k}  \\
& \le  & 2^{-0.47ht} 2^{t \log k}  \\
& \le & 2^{-0.46ht} \\
& < & \frac{1}{k(1.3)^{ht}}\;.
\end{eqnarray*}
As there are less than $k^s$ bijections $\pi$ with $s$ non-stationary points we obtain by the union bound,
the last two inequalities, and (\ref{e:pi-2}) that
$$
\Pr[H[K] \text{ and } H[J] \text{ are isomorphic} ] \le \sum_{s=2}^k 2^{-rk/2} k^s < \sum_{s=2}^k \frac{1}{k(1.3)^{ht}} <
\frac{1}{(1.3)^{ht}}\;.
$$
This completes the proof of Claim C2.
\qed

Let $H$ be a graph on $h$ vertices. The number of induced copies of $H$ in
a graph $G$ is denoted by $c_H(G)$. 
Let $c_H(n)$ denote the maximum of $c_H(G)$ taken over all graphs $G$ with $n$
vertices. It is easy to observe that for $h \le n \le 2h$ we have
$c_H(n) \ge 2^{n-h}$. Indeed, let $G$ be a graph obtained from $H$ by selecting
$n-h$ vertices of $H$ and duplicating them. Namely, selecting the vertices one by
one, if $v$ is a selected vertex, then add another vertex $v'$ and connect it
precisely to all the neighbors of $v$ in the current graph. This creates a graph
$G$ on $n$ vertices and a copy of $H$ in $G$ can be obtained by selecting
each non-duplicated vertex, and one of the two copies of each duplicated vertex.
The number of distinct copies of $H$ chosen in this way is $2^{n-h}$.
The following lemma shows that for graphs $H$ with $k(H) < \alpha h$ where
$\alpha < 1$ we in fact have $c_H(n) = 2^{n-h}$ at least when $n$ is not too large.

\begin{lemma}\label{l:num-copies}
Suppose $k(H) \le \alpha h$ where $0.5 < \alpha < 1$. Then for all $h \le n \le (2-\alpha)h$ we have $c_H(n) = 2^{n-h}$.
\end{lemma}
\Proof
Let $G$ be any graph on $n$ vertices. We prove that $c_H(G) \le 2^{n-h}$.
We will assume that the vertices of $H$ are labeled by $\{1,\ldots,h\}$.
The vertices of $G$ are labeled by $\{v_1,\ldots,v_n\}$.
We may associate a copy of $H$ in $G$ with an injection $f:[h]\rightarrow V(G)$.
Let ${\cal F}$ denote the set of all copies of $H$ in $G$.
Now, if $f,f' \in {\cal F}$, then $|Im(f)|=h$ and $|Im(f')|=h$,
thus $|Im(f) \cap Im(f')| \ge 2h-n \ge \alpha h$.
But since $k(H) \le \alpha h$ we have that if
$v \in Im(f) \cap Im(f')$, then $f^{-1}(v) = f'^{-1}(v)$.

Let $V^* = \cup_{f \in {\cal F}} Im(f)$. So, $V^* \subseteq V(G)$.
We may therefore assign a role to each $v \in V^*$, where the role of $v$ is $i$
if $f^{-1}(v)=i$ for some $f \in {\cal F}$. By the above paragraph,
roles are well-defined.

Hence $V^*$ may be partitioned into $V^*_1, \ldots, V^*_h$ where all the vertices
in $V^*_j$ have role
$j$. Now, every copy of $H$ in $G$ (i.e. every member of ${\cal F}$) is formed by
selecting one vertex from each $V^*_j$ for $j=1,\ldots,h$.
So, the number of copies of $H$ in $G$ is at most
$\pi_{j=1}^h |V^*_j| \le 2^{n-h}$.
\qed

Finally, we need the following simple lemma.
\begin{lemma}\label{l:imply}
Suppose that $K$ is a subgraph of $H$ on at least $k(H)$ vertices.
Then, $k(K) \le k(H)$.
\end{lemma}
\Proof
Every subgraph of $K$ on at least $k(H)$ vertices is also a subgraph of $H$
and hence is asymmetric. Any two subgraphs of $K$ on at least $k(H)$ vertices
are also two subgraphs of $H$ and hence are non-isomorphic.
\qed

Theorem \ref{t:1} follows immediately from the following lemma and from Lemma \ref{l:low-k}.
\begin{lemma}\label{l:1}
Let $\gamma > 0$ be a constant. For all $H$ sufficiently large, if $k(H) \le (1-\gamma) h$, then $H$ is avoidable.
\end{lemma}
\Proof
Applying Lemma \ref{l:dependence}, let $h-o(h) \le k \le h$ be such that $K_q$ has a decomposition into $q$ pairwise
edge-disjoint copies of $K_k$, where $q=k^2-k+1$.

Next, recall that $\C(H,k)$ is the set of all induced subgraphs of $H$ on $k$ vertices.
We will prove that $X=\C(k) \setminus \C(H,k)$ has the decomposition property for $q$.
Once we establish that, we are done since Lemma \ref{l:reduce} implies that $H$ is avoidable.

Hence, it remains to prove that any red-blue coloring of $K_q$ can be decomposed into edge disjoint copies of $K_k$
(the fact that it can is already stated in the first paragraph of this proof) but with the additional requirement
that in each copy of $K_k$ of this decomposition, the blue edges induce a subgraph which is not in $\C(H,k)$.

Suppose now that $K \in \C(H,k)$. First observe that since $|K|=k \ge h-o(h) \ge (1-\gamma) h \ge k(H)$, we have that $K$
is asymmetric. Furthermore, by Lemma \ref{l:imply}, 
$$
k(K) \le k(H) \le (1-\gamma) h \le (1-\frac{\gamma}{2})k\;.
$$
Using $\alpha = 1-\gamma/2$ in Lemma \ref{l:num-copies} we obtain that
for all $k \le n \le (1+\gamma/2)k$ we have $c_K(n) = 2^{n-k}$.

Let $Q$ be any graph on $q$ vertices (equivalently, a red-blue coloring of $K_q$).
We next prove that the density of $K$ in $Q$, namely $c_K(Q)/\binom{q}{k}$ satisfies
$$
\frac{c_K(Q)}{\binom{q}{k}} < \frac{1}{q\binom{h}{k}}\;.
$$
Assume the contrary. Then, for any $n$ such that $k \le n \le q$, we would have a subgraph $G$ of $Q$ on $n$ vertices
such that the density of $K$ in $G$ is at least $\frac{1}{q\binom{h}{k}}$, namely
$$
c_K(G) \ge \frac{\binom{n}{k}}{q\binom{h}{k}}\;.
$$
We shall use $n=\lfloor (1+\gamma/2) k \rfloor$. But then, $c_K(G) \le c_K(n)=2^{n-k}$.
To arrive at the desired contradiction we only need to show that
$$
2^{n-k} < \frac{\binom{n}{k}}{q\binom{h}{k}}\;.
$$
Now,
$$
\binom{n}{k} = \binom{\lfloor (1+\gamma/2)k \rfloor}{k} > \left( \frac{(1+\frac{\gamma}{2})^{1+\frac{\gamma}{2}}}
{(\frac{\gamma}{2})^{\frac{\gamma}{2}}} -o(1) \right)^k
$$
and
$$
q\binom{h}{k} \le q\binom{k+o(k)}{k} < (1+o(1))^k\;.
$$
Thus indeed,
$$
\frac{\binom{n}{k}}{q\binom{h}{k}} > \left( \frac{(1+\frac{\gamma}{2})^{1+\frac{\gamma}{2}}}
{(\frac{\gamma}{2})^{\frac{\gamma}{2}}} -o(1) \right)^k > 2^{\frac{\gamma}{2}k} \ge 2^{n-k}\;.
$$
Suppose the vertices of $Q$ are $\{1,\ldots,q\}$. Let ${\cal D}$ be some $K_k$-decomposition of $K_q$.
Hence $|{\cal D}|=q$ and any $R \in {\cal D}$ is an induced $k$-vertex subgraph of $Q$.
If each $R \in {\cal D}$ is an element of $\C(k) \setminus \C(H,k)$ we are done,
but the problem is that some $R$ might be isomorphic to some element of $\C(H,k)$.
For a permutation $\pi$ of $[q]$, let ${\cal D}_\pi$ be the $K_k$-decomposition of $K_q$
corresponding to the permutation. That is, each $R \in {\cal D}$ now corresponds to $R_\pi \in {\cal D}_\pi$
where $V(R_\pi)=\{\pi(v)~|~ v \in V(R)\}$. We will prove that there exists $\pi$ such that 
each $R_\pi \in {\cal D}_\pi$ is an element of $\C(k) \setminus \C(H,k)$.
As usual, it would be convenient to prove this counting argument using probabilistic language.

Suppose that $\pi$ is chosen uniformly among all permutations of $[q]$. For a fixed $K \in \C(H,k)$,
recall that the density of $K$ in $Q$ is less than $\frac{1}{q\binom{h}{k}}$.
As ${\cal D}$ (and thus ${\cal D}_\pi$) have $q$ elements, the probability that some element of ${\cal D}_\pi$
is isomorphic to $K$ is less than $\frac{1}{\binom{h}{k}}$. As there are at most $\binom{h}{k}$ elements
in $\C(H,k)$, we have that the expected number of elements of ${\cal D}_\pi$ that are isomorphic to some
element of $\C(H,k)$ is less than $1$. Hence, there exists $\pi$ such that 
each $R_\pi \in {\cal D}_\pi$ is an element of $\C(k) \setminus \C(H,k)$.
As $Q$ was an arbitrary graph of $q$ vertices, we have proved that $X=\C(k) \setminus \C(H,k)$ has the decomposition property for $q$.
\qed

\section{Proof of Theorem \ref{t:2}}

\subsection{A sufficient condition for the avoidability of ${\cal F}(S,k)$}

We now prove our main theorem of this section, from which Theorem \ref{t:2} can be obtained as a (nontrivial) corollary.
Recall that for a set of integers $S \subseteq \{0,\ldots,k-1\}$ we let ${\cal F}(S,k)$ be the set of all graphs on $k$ vertices whose degree set is
contained in $S$. The next theorem gives a sufficient condition for ${\cal F}(S,k)$ to be avoidable.
\begin{theorem}\label{t:suff}
Suppose that for every real parameter $x \in (0,1)$, the following linear system of three equations in the variables
$\{p_i ~|~ i \in \{0,\ldots,k-1\} \setminus S\}$ has a nonnegative solution\footnote{A nonnegative solution is a solution where each coordinate
is nonnegative.}:
\begin{eqnarray*}
\sum_{i \in \{0,\ldots,k-1\} \setminus S} \frac{x^{i-2}}{(i-2)!}\frac{(1-x)^{k-1-i}}{(k-1-i)!}p_i & = & 1\;, \\
\sum_{i \in \{0,\ldots,k-1\} \setminus S} \frac{x^i}{i!}\frac{(1-x)^{k-3-i}}{(k-3-i)!}p_i & = & 1\;, \\
\sum_{i \in \{0,\ldots,k-1\} \setminus S} \frac{x^{i-1}}{(i-1)!}\frac{(1-x)^{k-2-i}}{(k-2-i)!}p_i & = & 1\;.
\end{eqnarray*}
Then ${\cal F}(S,k)$ is avoidable.
\end{theorem}
The proof of Theorem \ref{t:suff} is based on the following lemma.
\begin{lemma}
\label{l:main}
Let ${\cal R}(S,k)$ be the complement of ${\cal F}(S,k)$, namely the set of all graphs on $k$ vertices whose degree set is
not contained in $S$. If the linear system of Theorem \ref{t:suff} has a nonnegative solution for every real parameter $x \in (0,1)$,
then for any graph $G$ we have $\nu^*_{{\cal R}(S,k)}(G) \ge \binom{n}{2}/\binom{k}{2}-o(n^2)$.
\end{lemma}
Notice that Lemma \ref{l:main} together with Lemma \ref{l:hr} immediately implies Theorem \ref{t:suff},
since we get that $\nu_{{\cal R}(S,k)}(G) \ge \binom{n}{2}/\binom{k}{2}-o(n^2)$ which means that ${\cal F}(S,k)$ is avoidable.

\vspace{3mm}
\noindent
{\bf Proof of Lemma \ref{l:main}}.
Let $k \ge 3$ and $S \subseteq \{0,\ldots,k-1\}$ be fixed.
As stated earlier, ${\cal F}(S,k)$ is the set of all graphs on $k$ vertices whose degree set is
contained in $S$ and ${\cal R}(S,k)$ is the complement of ${\cal F}(S,k)$.

Let $G$ be a graph with $V(G)=[n]$.
Our goal is to design a fractional packing $\phi$ from $\binom {G}{{\cal R}(S,k)}$ to $[0,1]$ such that
$|\phi| = \binom{n}{2}/\binom{k}{2}-o(n^2)$. This will prove that $\nu^*_{{\cal R}(S,k)}(G) \ge \binom{n}{2}/\binom{k}{2}-o(n^2)$
and yield a proof of Lemma \ref{l:main}.

We will construct $\phi$ as a sum of smaller fractional packings $\phi_v$ from $\binom {G}{{\cal R}(S,k)}$ to $[0,1]$,
one for each $v \in V(G)=[n]$. So,
$$
\phi= \sum_{v=1}^n \phi_v\;.
$$
We next define each $\phi_v$ and prove that $\phi$ satisfies the definition of a fractional packing, i.e. that (\ref{e:fp}) is
satisfied for each pair of distinct vertices $\{x,y\} \subset [n]$.

We first state a few properties that we require $\phi_v$ to have.
\begin{enumerate}
\item[P1.]
$\phi_v(H) > 0$ only if $H \in \binom {G}{{\cal R}(S,k)}$ and $v \in V(H)$.
\item[P2.]
$|\phi_v| = \frac{n-1}{k(k-1)}-o(n)$\;.
\item[P3.]
For any $x \in [n] \setminus v$, the sum of the values of $\phi_v(H)$ over all $H$ that contain the pair $\{x,v\}$ is $1/k-o_n(1)$.
In other words,
$$
\sum_{H \in \binom {G}{{\cal R}(S,k)} \,:\, \{v,x\} \subset V(H)} {\phi_v(H)} = \frac{1}{k}-o_n(1)\;.
$$
\item[P4.]
For any pair $\{x,y\} \subset [n] \setminus v$, the sum of the values of $\phi_v(H)$ over all $H$ that contain the pair $\{x,y\}$ is
$(k-2)/(k(n-2))-o(1/n)$. In other words,
$$
\sum_{H \in \binom {G}{{\cal R}(S,k)} \,:\, \{x,y\} \subset V(H)} {\phi_v(H)} = \frac{k-2}{k(n-2)}-
o\left(\frac{1}{n}\right)\;.
$$
\end{enumerate}
Let us see that if properties P2, P3, and P4 hold for each $\phi_v$ where $v \in [n]$, then indeed 
$|\phi| = \binom{n}{2}/\binom{k}{2}-o(n^2)$ and $\phi$ is a valid fractional packing.
First observe that by property P2, $|\phi|=n(\frac{n-1}{k(k-1)}-o(n))=\binom{n}{2}/\binom{k}{2}-o(n^2)$.
Next, consider some pair $\{x,y\} \subset [n]$. By P3, the sum of the values of $\phi_x$ over the elements that contain the pair
is $\frac{1}{k}-o_n(1)$. Likewise, the sum of the values of $\phi_y$ over the elements that contain the pair
is $\frac{1}{k}-o_n(1)$. By P4, for any $v \notin \{x,y\}$, the sum of the values of $\phi_v$ over the elements that contain the pair
is $(k-2)/(k(n-2))-o(1/n)$. So, the overall sum of values of $\phi$ over all elements that contain the pair is at most
$$
2\left(\frac{1}{k}-o_n(1)\right) + (n-2)\left( \frac{k-2}{k(n-2)}-
o\left(\frac{1}{n}\right)\right) = 1-o_n(1)\;.
$$
Hence, $\phi$ is a valid fractional packing with the claimed value.

We proceed to define $\phi_v$. Let us first set $\phi_v(H)=0$ for every $H \in {\cal R}(S,k)$ with $v \notin V(H)$.
This guarantees P1.
For each $i \in \{0,\ldots,k-1\} \setminus S$ let $q_{v,i}$ be a nonnegative real to be chosen later.
Now, consider any subset $W$ of $k-1$ vertices of $[n]\setminus \{v\}$. Clearly $W \cup \{v\}$ induces a subgraph
of $G$ on $k$ vertices which may or may not be in ${\cal R}(S,k)$.
Denote this subgraph by $G[v,W]$. If $G[v,W] \in {\cal R}(S,k)$ we must define $\phi_v(G[v,W])$.
Recall that $N(v)$ denotes the set of neighbors of $v$ in $G$.
Let $i=|W \cap N(v)|$ and clearly $0 \le i \le k-1$. Set
$$
\phi_v(G[v,W]) =
\begin{cases}
        0 & \text{if } i \in S\\
        q_{v,i} & \text{otherwise}\;.
\end{cases}
$$
Notice that we do need to consider the case $i \in S$ since it is possible that $i \in S$ while $G[v,W] \in {\cal R}(S,k)$.
 
We next define the values of the $q_{v,i}$. These values will depend on Properties P2,P3,P4, on $i$, and on $d(v)$, the degree of $v$ in $G$.
The number of elements $G[v,W]$ that received the weight $q_{v,i}$ is
the number of subsets $W$ of $k-1$ vertices of $[n]\setminus \{v\}$ such that $i=|W \cap N(v)|$,  which is
$$
\binom{d(v)}{i}\binom{n-1-d(v)}{k-1-i}\;.
$$
So, to satisfy P2 we must have
\begin{equation}\label{e:1}
\sum_{i \in \{0,\ldots,k-1\} \setminus S} \binom{d(v)}{i}\binom{n-1-d(v)}{k-1-i}q_{v,i} = \frac{n-1}{k(k-1)}-o(n)\;.
\end{equation}
Consider some edge $(v,x) \in E(G)$. How many elements $G[v,W]$ that contain the edge $(v,x)$ received the weight $q_{v,i}$?
For this to occur, $W$ must contain $i$ neighbors of $v$, while $x$ is one of those neighbors. Hence, the number of such elements is
$$
\binom{d(v)-1}{i-1}\binom{n-1-d(v)}{k-1-i}\;.
$$
To satisfy P3 we must therefore have that
\begin{equation}\label{e:2}
\sum_{i \in \{0,\ldots,k-1\} \setminus S} \binom{d(v)-1}{i-1}\binom{n-1-d(v)}{k-1-i}q_{v,i} = \frac{1}{k}-o_n(1)\;.
\end{equation}
Similarly, consider some non-edge $(v,x) \notin E(G)$. How many elements $G[v,W]$ that contain this non-edge received the weight $q_{v,i}$?
For this to occur, $W$ must contain $i$ neighbors of $v$, while $x \in W$ is not one of those neighbors. Hence, the number of such elements is
$$
\binom{d(v)}{i}\binom{n-2-d(v)}{k-2-i}\;.
$$
To satisfy P3 we must therefore have that
\begin{equation}\label{e:3}
\sum_{i \in \{0,\ldots,k-1\} \setminus S} \binom{d(v)}{i}\binom{n-2-d(v)}{k-2-i}q_{v,i} = \frac{1}{k}-o_n(1)\;.
\end{equation}
Consider some pair $\{x,y\} \subset [n] \setminus v$ such that both $(x,v) \in E(G)$ and $(y,v) \in E(G)$.
The number of elements $G[v,W]$ that contain this pair and received the weight $q_{v,i}$ is
$$
\binom{d(v)-2}{i-2}\binom{n-1-d(v)}{k-1-i}\;.
$$
To satisfy P4 we must therefore have that
\begin{equation}\label{e:4}
\sum_{i \in \{0,\ldots,k-1\} \setminus S} \binom{d(v)-2}{i-2}\binom{n-1-d(v)}{k-1-i}q_{v,i} = \frac{k-2}{k(n-2)}-o\left(\frac{1}{n}\right)\;.
\end{equation}
By similarly considering pairs $\{x,y\} \subset [n] \setminus v$ such that both $x,y$ are non-neighbors of $v$ we get that in order to
satisfy P4 we must have
\begin{equation}\label{e:5}
\sum_{i \in \{0,\ldots,k-1\} \setminus S} \binom{d(v)}{i}\binom{n-3-d(v)}{k-3-i}q_{v,i} = \frac{k-2}{k(n-2)}-o\left(\frac{1}{n}\right)\;.
\end{equation}
Finally, by considering pairs $\{x,y\} \subset [n] \setminus v$ such that exactly one of $x,y$ is a neighbor of $v$ we get that in order to
satisfy P4 we must have
\begin{equation}\label{e:6}
\sum_{i \in \{0,\ldots,k-1\} \setminus S} \binom{d(v)-1}{i-1}\binom{n-2-d(v)}{k-2-i}q_{v,i} = \frac{k-2}{k(n-2)}-o\left(\frac{1}{n}\right)\;.
\end{equation}
So, the question we remain with is whether we can find nonnegative reals $q_{v,i}$ such that equations (\ref{e:1}-\ref{e:6}) hold.
To simplify notation, let us set $x=d(v)/(n-1)$ and hence $(1-x)=(n-1-d(v))/(n-1)$. Also let $p_{v,i}=n^{k-2}q_{v,i}$.
Thus, in these terms, (\ref{e:1}-\ref{e:6}) become:
\begin{equation}\label{e:1p}
\sum_{i \in \{0,\ldots,k-1\} \setminus S} \frac{x^i}{i!}\frac{(1-x)^{k-1-i}}{(k-1-i)!}p_{v,i} = \frac{1}{k(k-1)}-o_n(1)\;. \tag{e1}
\end{equation}
\begin{equation}\label{e:2p}
\sum_{i \in \{0,\ldots,k-1\} \setminus S} \frac{x^{i-1}}{(i-1)!}\frac{(1-x)^{k-1-i}}{(k-1-i)!}p_{v,i} = \frac{1}{k}-o_n(1)\;. \tag{e2}
\end{equation}
\begin{equation}\label{e:3p}
\sum_{i \in \{0,\ldots,k-1\} \setminus S} \frac{x^i}{i!}\frac{(1-x)^{k-2-i}}{(k-2-i)!}p_{v,i} = \frac{1}{k}-o_n(1)\;. \tag{e3}
\end{equation}
\begin{equation}\label{e:4p}
\sum_{i \in \{0,\ldots,k-1\} \setminus S} \frac{x^{i-2}}{(i-2)!}\frac{(1-x)^{k-1-i}}{(k-1-i)!}p_{v,i} = \frac{k-2}{k}-o_n(1)\;. \tag{e4}
\end{equation}
\begin{equation}\label{e:5p}
\sum_{i \in \{0,\ldots,k-1\} \setminus S} \frac{x^i}{i!}\frac{(1-x)^{k-3-i}}{(k-3-i)!}p_{v,i} = \frac{k-2}{k}-o_n(1)\;. \tag{e5}
\end{equation}
\begin{equation}\label{e:6p}
\sum_{i \in \{0,\ldots,k-1\} \setminus S} \frac{x^{i-1}}{(i-1)!}\frac{(1-x)^{k-2-i}}{(k-2-i)!}p_{v,i} = \frac{k-2}{k}-o_n(1)\;. \tag{e6}
\end{equation}
It is not difficult to see that the six equalities (\ref{e:1p}-\ref{e:6p}) are linearly dependent and have rank at most $3$
even without the $o_n(1)$ allowed error term.
Indeed,
\begin{eqnarray}
(\ref{e:3p}) & = & \frac{k-1}{1-x} \cdot (\ref{e:1p}) - \frac{x}{1-x} \cdot (\ref{e:2p})\;, \label{e:f1} \\
(\ref{e:5p}) & = & \frac{(k-1)(k-2)}{(1-x)^2} \cdot (\ref{e:1p}) - \frac{2x(k-2)}{(1-x)^2} \cdot (\ref{e:2p}) + \frac{x^2}{(1-x)^2} \cdot (\ref{e:4p})\;, \label{e:f2}  \\ 
(\ref{e:6p}) & = & \frac{k-2}{1-x} \cdot (\ref{e:2p}) - \frac{x}{1-x} \cdot (\ref{e:4p}) \label{e:f3} 
\end{eqnarray}
So, (\ref{e:1p}), (\ref{e:2p}), (\ref{e:4p}) span the system of six equations.
It will be slightly more convenient to work with (\ref{e:4p}), (\ref{e:5p}), (\ref{e:6p}) as they all have the same right hand side.
They also span the six equations since (\ref{e:f3}) shows that (\ref{e:2p}) is spanned by (\ref{e:4p}), (\ref{e:6p}) 
and thus (\ref{e:f2}) shows that (\ref{e:1p}) is also spanned by (\ref{e:4p}), (\ref{e:5p}), (\ref{e:6p}) and thus (\ref{e:f1})
shows that (\ref{e:3p}) is spanned by (\ref{e:4p}), (\ref{e:5p}), (\ref{e:6p}) as well.

Finally, notice that the coefficients of the left hand side of each of (\ref{e:4p}), (\ref{e:5p}), (\ref{e:6p}) are exactly the
coefficients of the left hand sides of the equations stated in Theorem \ref{t:suff}. Since equations (\ref{e:4p}), (\ref{e:5p}), (\ref{e:6p})  have
the same right hand side, solvability is maintained if we normalize to require that each right hand side is $1$, as in the equations stated in Theorem \ref{t:suff}. Finally, as we have no control over $x$, and we require solvability for each $v \in V(G)$ (and different $v$'s may have
different degrees, thus different $x$'s) we need to ensure solvability for each $x \in (0,1)$. This prove Lemma \ref{l:main}.
\qed

\subsection{Sets that satisfy the conditions of Theorem \ref{t:suff}}

We start this section with an example showing that for some $S$, the linear system of Theorem \ref{t:suff} can only be
non-negatively solved for
all $x \in I \subset (0,1)$, where $I$ has positive measure strictly less than $1$. Hence, Theorem \ref{t:suff} cannot be applied to
such sets.

Consider the case $k=4$ and $S=\{2\}$. Observe that in this case, ${\cal F}(S,k)=\{C_4\}$.
The set of variables $\{p_i ~|~ i \in \{0,\ldots,k-1\} \setminus S\}$ is thus just $\{p_0,p_1,p_3\}$.
The system in Theorem \ref{t:suff} therefore becomes:
\begin{eqnarray*}
xp_3 & = & 1\;, \\
(1-x)p_0+ xp_1 & = & 1\;, \\
(1-x)p_1 & = & 1\;.
\end{eqnarray*}
This system has a nonnegative solution only if $x \in (0,\frac{1}{2}]$.

Table \ref{table:1} contains a list of all maximal sets\footnote{If Theorem \ref{t:suff} holds for a set $S$, then it clearly holds for any subset of $S$ as one can set any additional variables to zero.} $S$ for which Theorem \ref{t:suff} holds, for $5 \le k \le 11$.

\begin{table}
\centering
\begin{tabular}{|c|c|}
\hline
$q$ & $S$\\
\hline
$5$ & $\{2\}$\\ 
\hline
$6$ & $\{2\}\, \{3\}$\\ 
\hline
$7$ & $\{3\}\, \{2,4\}$\\ 
\hline
$8$ & $\{2,4\}\, \{2,5\}\, \{3,5\}$\\ 
\hline
$9$ & $\{3,5\}\, \{2,4,5\}\, \{2,4,6\}\, \{3,4,6\}$\\ 
\hline
$10$ & $\{2,4,6\}\, \{2,5,6\}\, \{3,4,6\}\, \{3,4,7\}\, \{3,5,6\}\, \{3,5,7\}\, \{2,4,5,7\}$\\ 
\hline
$11$ & $\{3,5,7\}\, \{2,4,5,7\}\, \{2,4,5,8\}\, \{2,4,6,7\}\, \{2,4,6,8\}\, \{2,5,6,8\}\, \{3,4,6,7\}\, \{3,4,6,8\}\, \{3,5,6,8\}$\\ 
\hline
\end{tabular}
\caption{All maximal sets $S$ that satisfy Theorem \ref{t:suff} for $5 \le k \le 11$.}
\label{table:1}
\end{table}

While the values in this table are verified by a computer program, one particular symmetric pattern that emerges is
$S=\{2,4,\ldots,k-3\}$ when $k$ is odd.
Our goal is to prove that this holds for all odd $k$, thereby proving Theorem \ref{t:2}.

\vspace{3mm}
\noindent
{\bf Proof of Theorem \ref{t:2}.}\,
To prove Theorem \ref{t:2} using Theorem \ref{t:suff}, we need to prove that for all $x \in (0,1)$, the system
\begin{eqnarray*}
\sum_{i \in \{0,1,3,\ldots,k-4,k-2,k-1\}} \frac{x^i}{i!}\frac{(1-x)^{k-3-i}}{(k-3-i)!}p_i & = & 1\;, \\
\sum_{i \in \{0,1,3,\ldots,k-4,k-2,k-1\}} \frac{x^{i-1}}{(i-1)!}\frac{(1-x)^{k-2-i}}{(k-2-i)!}p_i & = & 1\;, \\
\sum_{i \in \{0,1,3,\ldots,k-4,k-2,k-1\}} \frac{x^{i-2}}{(i-2)!}\frac{(1-x)^{k-1-i}}{(k-1-i)!}p_i & = & 1\;.
\end{eqnarray*}
has a nonnegative solution $(p_0,p_1,p_3,\ldots,p_{k-4},p_{k-2},p_{k-1})$.
Let us denote the matrix of coefficients by $A$ (so $A$ has $3$ rows and $(k+3)/2$ columns), and the vector of variables by $\hat{p}$.
So we need to prove that $A\hat{p}=J$ has a nonnegative solution where $J$ is the column vector $(1,1,1)$.
By the classical Farkas' Lemma \cite{farkas-1902} (or directly using linear programming duality), this holds
if and only if for any vector $y=(y_1,y_2,y_3) \in R^3$ such that $yA$ is nonnegative, we must have $y_1+y_2+y_3 \ge 0$.

So, suppose  that $yA$ is nonnegative. We must prove that $y_1+y_2+y_3 \ge 0$.
Consider first the product of $y$ with the first column of $A$. The first column of $A$ corresponds to $i=0$ so it is the column vector
$((1-x)^{k-3}/(k-3)!,0,0)$. As we assume that $yA$ is nonnegative, this implies that $y_1 \ge 0$.
Consider now the product of $y$ with the last column of $A$. The last column of $A$ corresponds to $i=k-1$ so it is the column vector
$(0,0,x^{k-3}/(k-3)!)$. As we assume that $yA$ is nonnegative, this implies that $y_3 \ge 0$.

We now consider the remaining $(k-1)/2$ inequalities of the form $yA_j \ge 0$ where $A_j$ is column $j$ of $A$ and $j=1,\ldots,(k-1)/2$.
We sum all of these $(k-1)/2$ inequalities. This sum is an inequality of the form $y_1f_1(x)+y_2f_2(x)+y_3f_3(x) \ge 0$.
Specifically,
\begin{eqnarray*}
f_1(x) = \sum_{i=1,3,\ldots,k-4} \frac{x^i}{i!}\frac{(1-x)^{k-3-i}}{(k-3-i)!}\;, \\
f_2(x) = \sum_{i=1,3,\ldots,k-2} \frac{x^{i-1}}{(i-1)!}\frac{(1-x)^{k-2-i}}{(k-2-i)!}\;, \\
f_3(x) = \sum_{i=3,5,\ldots,k-2} \frac{x^{i-2}}{(i-2)!}\frac{(1-x)^{k-1-i}}{(k-1-i)!}\;.
\end{eqnarray*}
Observe that $f_1(x)=f_3(x)$. So, we know that $(y_1+y_3)f_1(x)+y_2f_2(x) \ge 0$, that $y_1 \ge 0$ and that $y_3 \ge 0$.
This, in turn, implies that
$$
y_2 \ge -(y_1+y_3) \frac{f_1(x)}{f_2(x)}\;.
$$
Thus,
$$
y_1+y_2+y_3 \ge (y_1+y_3) \left[ 1 - \frac{f_1(x)}{f_2(x)} \right]\;.
$$
But observe that
$$
0 \le \frac{(1-2x)^{k-3}}{(k-3)!} = f_2(x)-f_1(x)
$$
so indeed $y_1+y_2+y_3 \ge 0$.
\qed

\section{Some unavoidable graphs}

We prove Proposition \ref{p:1}. The fact that $F={K_k}$ is unavoidable for every $k \ge 2$ is trivial.
We prove next that $F=\{K_{1,k-1}\}$ is unavoidable for each $k \ge 3$. Let $\alpha < (k^2/8+1) ^{-1}$ be a positive constant.
Consider a partition of $[n]$ into sets $A,B$ with $|A|=\alpha n$. Color $K_n$ by coloring all edges in $E(B)$ blue, and
all edges in $E(A) \cup E(A,B)$ red. Consider a $K_k$-packing of this $K_n$ which leaves $o(n^2)$ edges unpacked.
There are at most $|E(A)|=\binom{|A|}{2}$ elements in this packing that contain at least one edge of $E(A)$.
Any such element contains at most $k^2/4$ edges of $E(A,B)$. So, altogether, all of these elements contain at most
$\binom{|A|}{2}k^2/4$ edges of $E(A,B)$. But these do not cover all $|E(A,B)|=\alpha(1-\alpha)n^2-o(n^2)$ edges of $E(A,B)$ since
$$
\binom{|A|}{2}k^2/4 < \frac{\alpha^2}{2}n^2\frac{k^2}{4} < \alpha(1-\alpha)n^2 - o(n^2)\;.
$$
Hence, there is an element of the packing which contains no edge of $E(A)$ and does contain an edge of $E(A,B)$. This element is thus a
red $F=\{K_{1,k-1}\}$.

To see that $K_{2,3}$ is unavoidable, consider a partition of $[n]$ into sets $A,B$ with $|A|=n/2$. Color $K_n$ by coloring $E(A,B)$ red
and coloring $E(A) \cup E(B)$ blue. Consider a $K_5$-packing of this $K_n$ which leaves $o(n^2)$ edges unpacked.
Any element of the packing which contains an edge of $E(A,B)$ is either a red $K_{1,4}$ or a red $K_{2,3}$. They cannot all be red $K_{1,4}$
as otherwise, since any red $K_{1,4}$ occupies four red edges of $E(A,B)$ and at least $n^2/4-o(n^2)$ edges of $E(A,B)$ are packed,
there would have been $n^2/16-o(n^2)$ elements in the $K_5$-packing, but all together they would occupy $10n^2/16 - o(n^2)$ edges of $K_n$,
while the latter only has less than $n^2/2$ edges.
The same example shows that $K_{3,4}$ is unavoidable. Any element of the packing which contains an edge of $E(A,B)$ is either a red $K_{1,6}$ or a red
$K_{2,5}$ or a red $K_{3,4}$. They cannot all be red $K_{1,6}$ or red $K_{2,5}$ 
as otherwise, since any red $K_{1,6}$ or red $K_{2,5}$ occupies at most $10$ red edges of $E(A,B)$ and at least $n^2/4-o(n^2)$ edges of $E(A,B)$ are
packed, there would have been at least $n^2/40-o(n^2)$ elements in the $r$-packing, but all together they would occupy at least
$21n^2/40 - o(n^2)$ edges of $K_n$, while the latter only has less than $n^2/2$ edges.

The proof that $K_4^-$ is unavoidable is slightly more involved.
Consider a partition of $[n]$ into sets $A_1,A_2,A_3,A_4$ with $|A_i|=n/4$ for $i=1,2,3,4$.
Color $K_n$ by coloring $E(A_i)$ blue for $i=1,2,3,4$ and coloring $E(A_1,A_2)$ blue as well. All other edges are red.
Consider a $K_4$-packing ${\cal L}$ of this $K_n$ with leaves $o(n^2)$ edges unpacked.
We claim that ${\cal L}$ must contain a red $K_4^-$. Suppose it does not.
As each element in the packing consists of $6$ edges, we have that $|{\cal L}|=n^2/12-o(n^2)$.
We partition the elements of ${\cal L}$ into five types as follows.
Type $1$ elements have two vertices in $A_3$ and two vertices in $A_4$.
Type $2$ elements have three vertices in $A_3$ and one in $A_4$, or vice versa.
Type $3$ elements have all their four vertices in $A_3$ or all their four vertices in $A_4$.
Type $4$ elements have two vertices in $A_3$ and no vertex in $A_4$, or two vertices in $A_4$ and no vertex in $A_3$.
Type $5$ elements are all remaining elements.
Let $t_i$ be the number of elements of type $i$ for $i=1,\ldots,5$.

Consider a packed edge $(x,y)$ where $x \in A_3$ and $y \in A_4$ and the element $S \in {\cal L}$ containing $(x,y)$.
We claim that $S$ is entirely in $A_3 \cup A_4$. Indeed, otherwise, $S$ has at least one vertex in $A_1 \cup A_2$.
Suppose w.l.o.g. that it has a vertex in $A_1$. Then, no matter where the fourth vertex resides, we obtain a red $K_4^{-}$, a contradiction.
Thus, we have that for any packed edge $(x,y)$ where $x \in A_3$ and $y \in A_4$, the element of ${\cal L}$ containing it
is of type $1$ or of type $2$.
As there are $n^2/16-o(n^2)$ packed edges in $E(A_3,A_4)$, we have that $4t_1+3t_2=n^2/16-o(n^2)$.
Also, the number of edges of type $1$, type $2$ and type $3$ elements in $E(A_3) \cup E(A_4)$ is $2t_1+3t_2+6t_3$.

There remain $n^2/16-2t_1-3t_2-6t_3-o(n^2)$ edges in $E(A_3) \cup E(A_4)$ that are not of type $1,2,3$.
Hence, $t_4 \le n^2/16-2t_1-3t_2-6t_3$. Now, any element of type $5$ has at least three vertices in $A_1 \cup A_2$.
Also, each element of type $4$ contains a single edge of $E(A_1 \cup A_2)$.
As the number of edges in $E(A_1 \cup A_2)$ is less than $n^2/8$ we have that
$t_5 \le (n^2/8-t_4)/3$.
It follows that
\begin{eqnarray*}
|{\cal L}| & = & t_1+t_2+t_3+t_4+t_5\\
&  \le & t_1+t_2+t_3+t_4 + \frac{n^2/8-t_4}{3} \\
& = & t_1+t_2+t_3 +\frac{2}{3}t_4+\frac{n^2}{24}\\
& \le & t_1+t_2+t_3 +\frac{2}{3}\left( \frac{n^2}{16}-2t_1-3t_2-6t_3\right) + \frac{n^2}{24}\\
& = & -\frac{1}{3}t_1-t_2 -3t_3 + \frac{n^2}{12}\\
& = & t_1 - \frac{n^2}{48}+o(n^2)-3t_3 + \frac{n^2}{12}\\
& \le & t_1 - \frac{n^2}{48} + \frac{n^2}{12} +o(n^2)\\
& \le & \frac{n^2}{64} - \frac{n^2}{48} + \frac{n^2}{12} + o(n^2)
\end{eqnarray*}
contradicting the fact that $|{\cal L}|=n^2/12-o(n^2)$.
\qed

\section{Concluding remarks and open problems}

In the proof of Theorem \ref{t:1} we used Lemma \ref{l:low-k} that shows that $H \sim {\cal G}(h,\frac{1}{2})$ is highly
asymmetric, namely it has $k(H) \le \beta h$ for all $\beta \ge 0.94$, asymptotically almost surely.
However, it is not difficult to modify the proof of Lemma \ref{l:low-k}
so that it holds for $H \sim {\cal G}(h,p)$ for any constant $p \in (0,1)$. This would cause the lower bound for $\beta$ to increase towards $1$
(but staying strictly less than $1$),
changing some constants in the proof as the probability of the agreement event in the proof changes from $\frac{1}{2}$ to $p^2+(1-p)^2$.
Since in the proof of Lemma \ref{l:1} we can choose $\gamma$ to be any small positive constant, we obtain that
for every fixed $p \in (0,1)$, the random graph $H \sim {\cal G}(h,p)$ is avoidable asymptotically almost surely.

Theorem \ref{t:suff} gives a sufficient condition for avoidability of the family of graphs ${\cal F}(S,k)$,
namely all $k$-vertex graphs whose degrees are in $S$. It seems interesting to determine all maximal sets
$S \subset \{0,\ldots,k-1\}$ for which ${\cal F}(S,k)$ is avoidable. While this is trivial for $k=2,3$, the following
proposition determines the case $k=5$.
\begin{prop}\label{p:only}
$S=\{2\}$ and $S=\{1,3\}$ are the only maximal sets for which for which ${\cal F}(S,5)$ is avoidable.
\end{prop}
\Proof
The set ${\cal F}(\{1,3\},5)$ is trivially avoidable because it is empty (no graph with an odd number of vertices can have all its degrees odd).
The set ${\cal F}(\{2\},5)$ is avoidable by Theorem \ref{t:2}. The set ${\cal F}(\{2,3\},5)$ is unavoidable since it contains $K_{2,3}$
which is unavoidable by Proposition \ref{p:1}. Similarly, the complement of $K_{2,3}$ is unavoidable so ${\cal F}(\{1,2\},5)$ is unavoidable.
The sets ${\cal F}(\{0\},5)$ and ${\cal F}(\{4\},5)$ are unavoidable since $K_5$ and its complement are unavoidable.
Hence, $S=\{2\}$ and $S=\{1,3\}$ are the only maximal sets for which for which ${\cal F}(S,5)$ is avoidable.
\qed

Similar to the way Problem \ref{prob:1} asks to generalize the result of Erd\H{o}s and Hanani \cite{EH-1963},
it may be interesting to consider the analogous problem for exact decompositions, generalizing Wilson's Theorem.
Recall from Section 3 that $X \subseteq \C(k)$ has the decomposition property for $n$ if every red-blue coloring of $K_n$ has an
$X$-packing of size $\frac{n(n-1)}{k(k-1)}$. Accordingly, we say that $X \subseteq \C(k)$ has the {\em decomposition property}
if for all $n$ sufficiently large, $X$ has the decomposition property for $n$ whenever $\C(k)$ has the decomposition property for
$n$ (namely, by Wilson's Theorem, whenever $n \equiv 1,k \bmod k(k-1)$).
Similarly, we can define decomposition avoidability for graphs and sets. The following problem analogous to Problem 1 emerges.
\begin{prob}\label{prob:2}
For every fixed $k$, determine the subsets of $\C(k)$ that have the decomposition property.
\end{prob}
It is straightforward to see that if $H$ is decomposition avoidable, then it is also avoidable.
However, the following proposition might suggest that the converse is not true.
\begin{prop}\label{p:c4}
$C_4$ is decomposition unavoidable.
\end{prop}
\Proof
Let $n$ be such that $K_n$ has a $K_4$ decomposition (in fact, this is known to hold for {\em all} $n \equiv 1,4 \bmod 12$).
Partition the vertices of $K_n$ into two parts $A$ and $B$ of sizes $\lceil n/2 \rceil$ and $\lfloor n/2 \rfloor$.
Color $E(A,B)$ blue and all the other edges red. If our $K_4$-decomposition avoids a blue $C_4$, then any element
of this decomposition occupies at most $3$ blue edges. As there are $\lfloor n^2/4 \rfloor$ blue edges,
the decomposition must contain a least $\lfloor n^2/4 \rfloor/3$ elements. But this is impossible since it contains
precisely $n(n-1)/12$ elements.
\qed

\section*{Acknowledgment}
I thank the referee for valuable comments.


\begin{thebibliography}{10}
	
	\bibitem{BHP-2001}
	R.C. Baker, G.~Harman, and J.~Pintz.
	\newblock The difference between consecutive primes, ii.
	\newblock {\em Proceedings of the London Mathematical Society},
	83(03):532--562, 2001.
	
	\bibitem{EH-1963}
	P.~Erd{\H{o}}s and H.~Hanani.
	\newblock On a limit theorem in combinatorical analysis.
	\newblock {\em Publ. Math. Debrecen}, 10:10--13, 1963.
	
	\bibitem{ER-1963}
	P.~Erd{\H{o}}s and A.~R{\'e}nyi.
	\newblock Asymmetric graphs.
	\newblock {\em Acta Mathematica Hungarica}, 14(3-4):295--315, 1963.
	
	\bibitem{farkas-1902}
	J.~Farkas.
	\newblock Theorie der einfachen ungleichungen.
	\newblock {\em Journal f{\"u}r die reine und angewandte {M}athematik},
	124:1--27, 1902.
	
	\bibitem{GKLO-2016}
	S.~Glock, D.~K{\"u}hn, A.~Lo, and D.~Osthus.
	\newblock The existence of designs via iterative absorption.
	\newblock arXiv:1611.06827, 2016.
	
	\bibitem{goodman-1959}
	W.~Goodman.
	\newblock On sets of acquaintances and strangers at any party.
	\newblock {\em The American Mathematical Monthly}, 66(9):778--783, 1959.
	
	\bibitem{HR-2001}
	P.E. Haxell and V.~R\"odl.
	\newblock Integer and fractional packings in dense graphs.
	\newblock {\em Combinatorica}, 21(1):13--38, 2001.
	
	\bibitem{keevash-2014}
	P.~Keevash.
	\newblock The existence of designs.
	\newblock arXiv:1401.3665, 2014.
	
	\bibitem{KSV-2002}
	J.H. Kim, B.~Sudakov, and V.H. Vu.
	\newblock On the asymmetry of random regular graphs and random graphs.
	\newblock {\em Random Structures and Algorithms}, 21(3-4):216--224, 2002.
	
	\bibitem{rodl-1985}
	V.~R{\"o}dl.
	\newblock On a packing and covering problem.
	\newblock {\em European Journal of Combinatorics}, 6(1):69--78, 1985.
	
	\bibitem{szemeredi-1978}
	E.~Szemer\'edi.
	\newblock Regular partitions of graphs.
	\newblock In {\em Probl\`emes combinatoires et th\'eorie des graphes ({C}olloq.
		{I}nternat. {CNRS}, {U}niv. {O}rsay, {O}rsay, 1976)}, volume 260 of {\em
		Colloq. Internat. CNRS}, pages 399--401. CNRS, Paris, 1978.
	
	\bibitem{wilson-1975}
	R.M Wilson.
	\newblock Decomposition of complete graphs into subgraphs isomorphic to a given
	graph.
	\newblock {\em Congressus Numerantium}, 15:647--659, 1975.
	
	\bibitem{yuster-2005}
	R.~Yuster.
	\newblock Integer and fractional packing of families of graphs.
	\newblock {\em Random Structures and Algorithms}, 26(1-2):110--118, 2005.
	
\end{thebibliography}
\end{document}